# On the valuation of arithmetic–average Asian options: the Geman–Yor Laplace transform revisited

Peter Carr and Michael Schröder

**Abstract.** The Laplace transform approach of [**GY**] is a celebrated advance in valuing Asian options. Its insights are fundamental from both a mathematical and a financial perspective. In this paper, we discuss two observations regarding the financial relevance of its results. First, we show that the [**GY**] Laplace transform is not that of an Asian option price, as reported in [**GY**] and other papers. We nonetheless show how the [**GY**] Laplace transform can be used to obtain the price of an Asian option. Second, we find that following [**GY**] these Laplace transfoms are available only if the risk-neutral drift is not less than half the squared volatility. Using complex analytic techniques, we lift this restriction, thus extending the financial applicability of the Laplace transform approach.

## 1. Introduction

Asian options are widely used financial derivatives whose valuation has intrigued finance theorists for over a decade. In general, these options provide nonlinear payoffs on the arithmetic average of the price of an underlying asset. Our understanding of Asian options has evolved as an intricate interplay between theoretical and computational perspectives. As early as 1992, Yor [**Y**] expressed the Black–Scholes prices of Asian options as a certain triple integral. While this remains one cornerstone of our theoretical understanding of Asian options, the Laplace transform approach of [**GY**] had far reaching consequences for the way we see this valuation problem today. It has first spurred the developement of an impressive array of increasingly sophisticated numerical inversion techniques, see [**FMW**] for a recent example. It has also motivated the development of analytic inversion techniques. For example, by analytically inverting the Laplace transform, [**SA**] improves on Yor's result by expressing the Black–Scholes price of an Asian option as a single integral of the product of two functions. To compute this integral, [**SE**] has derived series and asymptotic expansions. The other principal computational alternatives available for valuing Asian options are continued fractions, expansion into eigenfunctions, and expansion into orthogonal polynomials. Based on Yor's results in [**Y**], this last alternative has been pursued in [**Du**]. This work has made it possible to establish surprising theoretical links between Asian options and modular forms, see [**Y3**], [**SL**]. Unfortunately, on the practical side [**SL**] found computing with these Laguerre series to be prone to serious structural difficulties. Thus, at present, all workable non–PDE methods for computing exact Black–Scholes prices of Asian options derive from the Laplace transform of [**GY**].

It is widely believed that this Laplace transform is that of an Asian option, as can be seen from the original paper [**GY**] or [**FMW**] for example. This paper shows that this is not the case. Nonetheless, the [**GY**] Laplace transform can be used to obtain the price of an Asian option. A more detailed discussion of these somewhat paradoxical findings is deferred to Part I. To paraphrase those findings here, the function which is Laplace transformed in [**GY**] is a normalized option value in two variables. Both variables depend on time to maturity, so the variables are not independent and in fact the second depends on the first. This dependency is not accounted for when in [**GY**] the transform is taken over the first variable.

Thus, the first achievement of the paper is that we show how the [**GY**] Laplace transform can nonetheless be used to value Asian options. The key idea is to replace the problem of valuing the Asian option with the problem of valuing a whole family of auxiliary non–Asian options. We can then reconstruct the Black–Scholes price of the Asian option using this family as a whole.

Encouraged by the kind support of Yor, we have extended these results from the zero dividend yield situation originally considered in [**GY**] to the case of constant positive dividend yields. Part II gives an exposition of this extension, incorporating Yor's tutorials and suggestions. We count it as the second contribution of this paper to the understanding of Asian options.

In providing this extension, we found that the probabilistic arguments of [**GY**] are sufficient to establish the Laplace transforms only if the risk-neutral drift is not less than half the squared volatility. This is due to the still limited knowledge of Bessel processes with negative indices. Unfortunately, this lower bound on the drift restricts the financial applicability of the results. For example, if volatility is 30%, we simply do not know the Laplace transforms for the non–Asian options if the difference between the riskfree rate and the dividend yield is less than 4.5%. Worse yet, the greater is the volatility, the greater is this range of parameters in which we do not have the Laplace transforms. Unfortunately, it is precisely due to high volatility that Asian options are used in the first place.

Thus, our third contribution is to lift these restrictions in showing the existence of the desired Laplace transforms in Part III. As described in §7, we are thus now able to value Asian options for all constellations of the riskfree rate, dividend yield, and volatility. To establish our result we change perspective form probabilistic methods to complex analysis. We apply analytic continuation which is is a standard method from complex analysis. However, the analytic continuation of Laplace transforms in the present setting appears to be a highly non–standard way of applying this technique, and it has to be handled with caution and care. An overview of the key steps of this extension is in §7, while the pertinent background from complex analysis is collected in the appendices. The reader who perseveres with this article will hopefully join in our realization that valuing the Asian option is akin to valuing non–Asian options on the edge of a knife.

# Part I

## 2. Basic notions

First we discuss the notions basic for the analysis of Asian options. We work in the Black–Scholes framework using the risk–neutral approach to the valuation of contingent claims. In this set–up there are two securities. There is a riskless security, a bond, whose price grows at the continuously compounding positive interest rate $r$. There is also a risky security, whose price process $S$ is modelled as follows. Consider a complete probability space equipped with the standard filtration of a standard Brownian motion on the time set $[0, \infty)$. On this filtered space, we have the risk neutral measure $Q$, which is a probability measure equivalent to the given one. And then we have a standard $Q$–Brownian motion $B$ such that $S$ is the strong solution of the following stochastic differential equation:

$$dS_t = \varpi \cdot S_t \cdot dt + \sigma \cdot S_t \cdot dB_t, \qquad t \in [0, \infty).$$

The positive constant $\sigma$ is the volatility of $S$. The specific form of the otherwise arbitrary constant $\varpi$ depends on the nature of the security modelled (eg. stock, currency, commodity etc.). For example it is the interest rate if $S$ is a non–dividend–paying stock.

Fix any time $t_0$ and consider the accumulation process $J$ given for any time $t$ by:

$$J(t) = \int_{t_0}^{t} S_u \, du \, .$$

The European–style *arithmetic–average Asian option* written at time $t_0$, with maturity $T$, and fixed–strike price $K$ is then the contingent claim on the closed time interval from $t_0$ to $T$ paying

$$\left( \frac{J(T)}{T - t_0} - K \right)^+ = \max \left\{ 0, \frac{J(T)}{T - t_0} - K \right\}$$

at time $T$. Recall that points in time are taken to be non–negative real numbers.

## 3. Normalized valuation of the Asian option

In the setting of the previous section, the price $C_t$ of the Asian option at any time $t$ between $t_0$ and $T$ is given as the following risk neutral expectation

$$C_t = e^{-r(T-t)} E^Q \left[ \left( \frac{J(T)}{T - t_0} - K \right)^+ \middle| \mathscr{F}_t \right]$$

which is conditional on the information $\mathscr{F}_t$ available at time $t$. However, following [**GY**, §3.2], do not focus on this price. As described there in great detail, we instead normalize the valuation problem as follows. Multiply and divide the strike by the length of the monitoring period, then factor out the reciprocal of this length. Split the integral $J(T)$ into two integrals, one of which is deterministic by time $t$ and the other of which is random. Couple the deterministic integral with the new strike. For the random integral, restart the Brownian motion driving the underlying at time $t$, and then using the scaling property of Brownian motion, change time to normalize the volatility in the new time scale to twohundred percent. The precise result is the factorization:

$$C_t = \frac{e^{-r(T-t)}}{T-t_0} \cdot \frac{4S_t}{\sigma^2} \cdot C^{(\nu)}(h,q)$$

which reduces the general valuation problem to computing

$$C^{(\nu)}(h,q) = E^Q\big[(A_h^{(\nu)} - q)\big],$$

the *normalized time–t price* of the Asian option. To explain the notation, $A^{(\nu)}$ is Yor's twohundred percent volatility accumulation process

$$A_h^{(\nu)} = \int_0^h e^{2(B_w + \nu w)} dw \,,$$

and the normalized parameters are as follows:

$$\nu = \frac{2\varpi}{\sigma^2} - 1, \qquad h = \frac{\sigma^2}{4}(T-t), \qquad q = kh + q^*,$$

where

$$k = \frac{K}{S_t}, \qquad q^* = q^*(t) = \frac{\sigma^2}{4S_t}\left(K \cdot (t-t_0) - \int_{t_0}^t S_u\,du\right).$$

To interpret these quantities, $\nu$ is the *normalized adjusted interest rate*, $h$ is the *normalized time to maturity*, which is non–negative, and $q$ is the *normalized strike price*. Notice how $q$ becomes affine linear in the time variable $h$ with coefficients $k$ and $q^*$ depending only on quantities known at time–$t$.

## 4. The pricing dichotomy

There is now a dichotomy in computing the normalized time–$t$ price

$$C^{(\nu)} := E^Q\big[(A_h^{(\nu)} - q)^+\big]$$

of the Asian option according to the normalized strike price $q$ being positive or not. Indeed, if $q$ is non–positive, Asian options loose their option feature. Computing their normalized time–$t$ prices $C^{(\nu)}$ is straightforward:

$$C^{(\nu)} = E^Q\Big[A_h^{(\nu)}\Big] - q \,.$$

On applying Fubini's theorem, this last expectation is computed as follows:

$$E^Q\Big[A_h^{(\nu)}\Big] = \frac{e^{2h(\nu+1)} - 1}{2(\nu+1)},$$

for $\nu$ any real number. The right hand side is analytic in $\nu$ with its value at $\nu = -1$ being equal to $h$.

In this way, computing $C^{(\nu)}$ is reduced to the case where $q$ is positive. As a first step in this direction, the aim of the present paper is to develop a certain family of Laplace transforms. The valuation problem then reduces to one of inversion. This inversion has been effected analytically in [**SA**] with an integral for the normalized price as the main result. Series and asymptotic expansions for computing this integral are given in [**SE**].

## 5. Intuition about the Laplace transform

The first key idea for computing the normalized time–$t$ price

$$C^{(\nu)} = E^Q \Big[ \big( A_h^{(\nu)} - (kh + q^*) \big)^+ \Big]$$

of the Asian option is, roughly speaking, to make the maturity date a variable. This is the starting point for the Laplace transform approach of [**GY**]. Since the mathematical background is summarized in §17, we try to provide financial intuition here.

Making the maturity a variable translates into considering $C^{(\nu)}$ via $h(y) = (\sigma^2/4) \cdot (y - t)$ as a function of the maturity $y$ ranging from $t$ to $\infty$. Since option values increase with maturity, it seems reasonable to give lower weight to more distant maturities. Translate this intuition by continuously discounting $y \mapsto C^{(\nu)}(h(y))$ with a constant discount rate $z$. There is of course no natural choice for $z$, so we consider it as a variable. In this way, we associate to the function of maturity $y \mapsto C^{(\nu)}(h(y))$ a function of this discount rate $z$. Written in full length, this *Laplace transform* of $C^{(\nu)}$ is thus formally given by

$$\mathscr{L}(C^{(\nu)})(z) = \int_0^\infty e^{-zx} E^Q \Big[ \big( A_x^{(\nu)} - (kx + q^*) \big)^+ \Big] dx$$

on changing variables $x = h(y) = (\sigma^2/4) \cdot (y - t)$.

It is a characteristic feature of the "normalized strike price" of the Asian option that it depends on the transformation variable $x$ via $q(x) = kx + q^*$. In [**GY**], however, options are considered whose "normalized strike price" is independent of $x$, and, as discussed in §9, this is crucial for the validity of their Laplace transform computations. So it appears initially that their results cannot be used to value Asian options. Luckily, they can nevertheless be used for exactly this purpose, and this is explained next.

## 6. Valuing Asian options using families of non–Asian options

The basic idea is as follows. Try to reconstruct the normalized time–$t$ price

$$C^{(\nu)} = E^Q \Big[ \big( A_h^{(\nu)} - (kh + q^*) \big)^+ \Big]$$

of the Asian option from a family of auxiliary functions whose single members are unrelated to the problem of valuing the Asian option but amenable to the [**GY**] analysis. The picture for this is to think of the graph of the normalized price being fibred up: try to put through any point of this graph a curve that is amenable to the analysis of [**GY**]. Given that the time dependency of the normalized strike price is the obstacle to applying these results, there are natural candidates for such functions. Force this strike price to be constant, and consider for any real number $a$ the function $f_{GY,a}$ on the positive real line that sends any $x > 0$ to

$$f_{GY,a}(x) = E^Q\big[(A_x^{(\nu)} - a)^+\big].$$

They are the prices of certain non–Asian options. Taken individually, they cannot be used to value the original Asian option. However, as a whole they allow one to recover the normalized time–$t$ price. In the notation of §3 this is made precise in the following

**Key Reduction:**   *If $q = kh + q^*$ is positive, computing the normalized time–$t$ price $C^{(\nu)}$ of the Asian option reduces to computing all $f_{GY,a}$ with $a > 0$. More precisely, $C^{(\nu)}$ is obtained by choosing the function $f_{GY,kh+q^*}$ and evaluating it at $h$.*

A moment's reflection will convince the reader that this is true by construction. In conclusion we would like to stress that the basic idea of our valuation technique is to value an option using a whole family of auxiliary options. And of course, this works in much more general situations.

# 7. The Laplace transform of $f_{GY,a}$

The key reduction of the preceding §6 shows the way to apply the [**GY**] Laplace transform to value Asian options. Moreover, adopting the notation of §3, consider the family of all functions $f_{GY,a}$ with $a > 0$ that send any $x > 0$ to

$$f_{GY,a}(x) = E^Q\big[(A_x^{(\nu)} - a)^+\big].$$

Recall that its single members are unrelated to valuing the Asian option, but as a whole, the family allows reconstruction of the normalized time–$t$ price $C^{(\nu)}$ of the Asian option. As first step in this direction, try to compute the Laplace transform of any $f_{GY,a}$

$$F_{GY,a}(z) = \int_0^\infty e^{-zx} f_{GY,a}(x)\,dx = \mathscr{L}(f_{GY,a})(z).$$

Here the complex number $z$ is to be taken in a half–plane sufficiently deep within the right complex half–plane such that the integrals are finite. The function so obtained is analytic. The precise conditions under which these integrals are finite is part of the description of these *generalized modified Geman–Yor Laplace transforms* in our

**Theorem:**   *Setting $\nu = 2\sigma^{-2}\varpi - 1$, the integrals $F_{GY,a}$ are finite for any complex number $z$ with $\mathrm{Re}\,(z) > \max\{0, 2(\nu+1)\}$. Let $z$ be any complex number with positive real part $z_0$ such that $z_0 > 2(\nu+1)$ if $\nu \geq 0$ and $z_0 > 4$ if $\nu < 0$. Then we have*

$$F_{GY,a}(z) = \frac{D_\nu(a, z)}{z(z - 2(\nu+1))}$$

*where on choosing the principal branch of the logarithm*

$$D_\nu(a, z) = \frac{e^{-\frac{1}{2a}}}{a} \int_0^\infty e^{-\frac{x^2}{2a}} x^{\nu+3} I_{\sqrt{2z+\nu^2}}\Big(\frac{x}{a}\Big)\,dx\,.$$

Herein $I_\mu$ is the modified Bessel function with complex order $\mu$, as discussed in [**L**, Chapter 5]. The integral $D_\nu(a, z)$ can be expressed using the confluent hypergeometric function $\Phi$ discussed in [**L**, Chapter 9]. The precise result is the following

**Corollary:** *Let $z$ be any complex number with positive real part $z_0$ such that $z_0 > 2(\nu+1)$ if $\nu \geq 0$ and $z_0 > 4$ if $\nu < 0$. Then we have*

$$F_{GY,a}(z) = \frac{D_\nu(a, z)}{z(z - 2(\nu+1))}$$

*where*

$$D_\nu(a, z) = \frac{\Gamma\left(\dfrac{\nu+4+\mu(z)}{2}\right)}{\Gamma\left(\mu(z)+1\right)} \cdot \Phi\left(\frac{\nu+4+\mu(z)}{2}, \mu(z)+1; \frac{1}{2a}\right) \cdot \frac{(2a)^{\frac{\nu+2-\mu(z)}{2}}}{e^{\frac{1}{2a}}}$$

*on setting $\mu(z) = \sqrt{2z + \nu^2}$.*

We again stress that these Laplace transforms are not those of the Asian option price, but rather the Laplace transforms of the prices of auxiliary options. To obtain Asian option prices, we have to invert these Laplace transforms as indicated in §17. In full mathematical generality, analytic inversion has been achieved in [**SA**]. Numerical inversions have also been accomplished. For example [**FMW**] computed the following seven cases as reproduced in [**Du**, Table 7.1]

| Case | $r$ | $\sigma$ | $T$ | $S_0$ | $2C^{(\nu)}$ |
|------|------|------|-----|-------|--------------|
| 1 | 2% | 10% | 1 | 2.0 | 0.056 |
| 2 | 18% | 30% | 1 | 2.0 | 0.219 |
| 3 | 1.25% | 25% | 2 | 2.0 | 0.172 |
| 4 | 5% | 50% | 1 | 1.9 | 0.194 |
| 5 | 5% | 50% | 1 | 2.0 | 0.247 |
| 6 | 5% | 50% | 1 | 2.1 | 0.307 |
| 7 | 5% | 50% | 2 | 2.0 | 0.352 |

**Table 1** Prices $2C^{(\nu)}$ for $K = 2.0$ and
$t_0 = t = 0$ using numerical Laplace inversion

However, the results of [**GY**] have been applied literally in computing these values. So it is only our two results that give a rigorous base for these computations and justify their principal relevance and validity for valuing the Asian option *a posteriori*.

Obtaining our two results proceeds in two steps. In a first probabilistic step, the arguments of [**GY**] are adapted to compute the modified Geman–Yor transforms $F_{GY,a}$. A new exposition of their argument is in Part II incorporating the tutorials and kind suggestions of Yor. Its key idea is to factorize the geometric Brownian motion of the underlying over a Bessel process of index $\nu$. Pertinent notions are discussed in §8. This makes time stochastic in such a way that Yor's accumulation process now takes the double role of both a stochastic clock and a control variable for the Asian option. At first sight, this appears to complicate the original valuation problem. However, this double role of Yor's accumulation process is especially suited to the Laplace transform. Indeed, using the fact that the "strike price" $a$ is independent of time makes it possible to reduce the computation of the Laplace transform to the following problem: obtain explicit expressions for the Bessel semigroup of index $\nu$ and for a certain conditional expectation involving first passage times of Yor's accumulation process.

However, such results are available for both concepts only if the index $\nu$ is non–negative. This non–negativity condition translates into the condition that the risk-neutral drift is not less than half the squared volatility. If this condition is violated, we do not have the Laplace transform. Unfortunately, this places restrictions on the Finance applicability of the result as dicussed in the introduction. Their effect can already be seen from the fact that $\nu$ is non–negative only in the first two numerical inversions reproduced in Table 1 above.

In a second analytic step, we thus lift this non–negativity restriction on $\nu$ in §10 to §15. For this we use analytic continuation which is a technique from complex analysis. To get an idea how it works, the Theorem and the Corollary both ask for establishing the equality of two functions. The trick is to treat $\nu$ as a variable too and to consider these as functions of both the transformation variable $z$ and $\nu$. In fact, the trick is to forget the transformation variable as much as possible and consider them as functions in $\nu$ alone. From that point of view first part of the argument establishes the equality of two functions in $\nu$ on the non–negative real line.

The problem is to establish an analogous statement for $\nu$ varying over the negative real line, and this is what we do in Part III. The idea is to prove more than that. This is because functions become amenable to complex analytic methods only on open subsets of the complex plane. Thus the *identity theorem* of §16.2 which we want to apply asserts that two functions on a connected open subset of the complex plane are equal if they are analytic and agree on a convergent sequence only. So it is in fact no longer possible to stick to real numbers only. As a subset of the complex plane they are closed with an empty interior, and we have to find ways to thicken the real line. At this stage, however, nothing is known about the functions if $\nu$ is outside the non–negative real line. In particular, it is not known if they exist at all, and this we have to establish together with their analyticity properties.

A number of principal problems come up here that have to be be handled with care and caution. For instance, the domains of definition in $z$ of these Laplace transforms change with the continuation variable $\nu$. We have been able to cope with that only by establishing yet another upper bound on the price of the Asia option.

The actual extension argument of Part III is by a two–front attack. The evaluation of $D_\nu(a, z)$ as a generalized Weber integral using the confluent hypergeometric function $\Phi$ is classical if $\nu$ is non–negative and $z$ is real. The point is that the identity so obtained remains true for complex numbers $z$ and $\nu$ in the following situation. The real part of $z$ is positive and $\nu$ varies over a band of a certain height depending on $\mathrm{Re}\,(z)$ which is symmetric with respect to the real axis. On this domain, the functions involved are analytic in $\nu$ as explained in §11. Crucial for this is our somewhat delicate square root lemma of §10. And at this stage we have the Corollary given the Theorem is valid.

The analyticity properties in $\nu$ of the risk neutral expectation $L(x, \nu) = E^Q[(A_x^{(\nu)} - a)^+]$. are at the base of the Theorem. Originally defined only for non–negative indices $\nu$ we show that this function extends to the whole complex plane as an analytic function in $\nu$. This determines the analyticity properties of the Laplace transform $F(\nu)(z) = \mathscr{L}(L(\cdot, \nu))(z)$ which are discussed in §14. The key idea is to consider for any $\varepsilon > 0$ the strip in the plane of height $2\varepsilon$ symmetric with respect to the real axis, the $\varepsilon$–thickened real line. On choosing $z$ with $\mathrm{Re}\,(z)$ bigger than four and $\varepsilon$ sufficiently small we can show analyticity of $\nu \mapsto F(\nu)(z)$ for elements $\nu$ in this $\varepsilon$–thickened real line with real part smaller than $\varepsilon$. This $\varepsilon$–thickened half–line contains the open interval $(0, \varepsilon)$. This makes possible the analytic continuation in $\nu$ of the Laplace transforms. Indeed, as discussed in §15, having the Laplace transform identity of the Theorem for $\nu$ in this open interval is sufficient to extend its validity to the $\varepsilon$–thickened half–line and so to all negative indices $\nu$. Drawing a picture of the domains, valuing the Asian option so is akin to valuing non–Asian options on the edge of a knife.

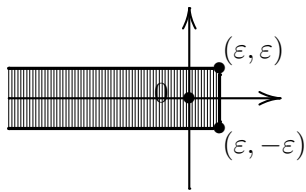

**Figure 1** The $\varepsilon$–thickened half–line

# Part II  Probabilistic part of the proof

## 8. Preliminaries on Bessel processes

As a preliminary to establishing the Laplace transforms of the preceding section, this section collects a number of pertinent facts from the theory of Bessel processes. They are patterned after the example of the Bessel processes of integer dimension $\delta \geq 2$ which are defined by taking the Euclidian distance from the origin of a Brownian motion in dimension $\delta$. Applying the Itô Lemma, their squares are the continuous strong solutions of the stochastic differential equation

$$d\rho_t = 2\delta dt + 2\sqrt{\rho_t}\, dB_t, \qquad \rho_0 = 1$$

[**RY**, XI, §1]. These stochastic differential equations make sense for any real number $\delta$ and have a unique continuous strong solution also if $\delta$ is smaller than two. The latter is by definition the square of the *Bessel process* $R^{(\nu)}$ on $[0, \infty)$ with *index* $\nu = \delta/2 - 1$ starting at 0 with the value 1. While Bessel processes of non–negative indices $\nu$ stay positive if started with a positive value at time zero, Bessel processes of negative indices $\nu$ develop some pathologies. As explained in [**RY**, XI, §1], in this case they hit zero. If $-1 < \nu < 0$, they are thereupon instantenously reflected and never become negative.

For $\nu = -1$ they continue zero. If $\nu < -1$, they even continue negative. This is why the following *Lamperti identity* may be surprising

**Lemma:** *For the index $\nu$ any real number, we have the factorization:*

$$e^{B_t + \nu t} = R^{(\nu)}\big(A_t^{(\nu)}\big).$$

*for any $t > 0$, where $A_t^{(\nu)} = \int_0^t e^{2(B_w + \nu w)} dw$ is Yor's accumulation process.*

For $\nu \geq 0$ a proof is given in [**Y92**, §2] while the general case is now contained as exercise XI (1.28), p.452 in the third edition of [**RY**]. We are indebted to Yor for this and for kindly supplying us with the following argument.

The trick is to apply the Itô rule to the square $Z_t$ of $Y_t = \exp(\nu t + B_t)$ to get

$$Z_t = 2(\nu + 1) \int_0^t Z_w dw + 2 \int_0^t Z_w dB_w \,.$$

Time change the process using the inverse function $\tau(t) = \inf\{u | \int_0^u Z_w \, dw > t\}$ to Yor's accumulation process to get

$$Z_{\tau(t)} = 2(\nu + 1)t + 2 \int_0^{\tau(t)} Z_w dB_w \,.$$

To interpret the stochastic integral in this sum, apply the basic time change formalism for stochastic processes as in [**Ø**, §8.5] to obtain

$$\int_0^{\tau(t)} Z_w dB_w = \int_0^t Z_w \sqrt{\tau'(w)} \, dW_w$$

where $W_t$ is defined as the stochastic integral $W_t = \int_0^{\tau(t)} \sqrt{Z_w} \, dB_w$ and is a Brownian motion. Using the inverse function theorem of calculus, the derivative of $\tau$ is equal to the reciprocal of the derivative with respect to time of Yor's accumulation process at time $w$. Hence it is equal to the reciprocal of $Z_w$. On substitution we so identify the time changed process $Z$ as continuous solution to the stochastic differential equation for the square of the Bessel process of index $\nu$:

$$Z_{\tau(w)} = 2(\nu + 1)t + 2 \int_0^t \sqrt{Z_w} \, dW_w \,.$$

Using the uniqueness of the solution of these stochastic differential equations, the time–changed process $Z$ so is the square of the Bessel process of index $\nu$. Reversing the time change, this translates into

$$Y_t^2 = (R_t^{(\nu)})^2 (A_t^{(\nu)}).$$

To establish the identity of the Lemma we have to take square roots. This is not a problem if $\nu$ is non–negative since then the Bessel process takes non–negative values only. It does pose a problem if $\nu$ is negative. In this case, however, recall that the Bessel process starts at time zero with a positive value. Since it is continuous by hypothesis, it will stay positive until it first hits zero at time $t^* > 0$. Since the accumulation process starts with zero at time zero there is a latest point in time $t^{**}$, infinity admitted, such that $A^{(\nu)}$ is smaller than $t^*$ at all points in time $t$ smaller than $t^{**}$. Thus we have the required identity at least for all points in time $t$ smaller than $t^{**}$. Now $Y_t$ is never zero. Since the processes on both sides of the identity are continuous in time, $t^{**}$ must be infinity, and the proof is complete.

## 9. The modified Geman–Yor Laplace transform for $\nu \geq 0$

This section is the first step in the proof of the integral representation of §7 Theorem for the Laplace transform

$$F_{GY,a}(z) = \int_0^\infty e^{-zx} f_{GY,a}(x)\,dx = \mathscr{L}(f_{GY,a})(z),$$

where in the notation of §3 and §6 we have $f_{GY,a}(x) = E^Q[(A_x^{(\nu)} - a)^+]$, for any positive real numbers $a$ and $x$. We explain how under the restriction $\nu = 2\sigma^{-2}\varpi - 1 \geq 0$ the probabilistic arguments of [**GY**] apply *mutatis mutandis* to give the

**Lemma:** *The assertions of §7 Theorem are valid if $\nu = 2\sigma^{-2}\varpi - 1 \geq 0$.*

We are very indebted to Yor for correspondence and discussions about this result, and are very grateful for his kind support. In the sequel we want to indicate the key steps of the proof following [**GY**] but trying to incorporate his tutorials and suggestions. Hopefully no pitfalls have remained undetected.

The basic idea is to make time stochastic using the Lamperti identity

$$e^{\nu w + B_w} = R^{(\nu)}\big(A_w^{(\nu)}\big)$$

for all positive real numbers $w$ as it has been discussed in the preceding section. Here $R^{(\nu)}$ is the Bessel process of index $\nu$ with $R^{(\nu)}(0) = 1$. On applying this Lamperti identity, $A^{(\nu)}$ has the double role of both control variable and stochastic clock. That the "strike price" $a$ is independent of time now becomes essential. It makes possible to transcribe the condition on the control variable to be bigger than $a$ as the first passage time

$$\tau_{\nu,a} = (\text{first passage time of } A^{(\nu)} \text{ to level } a)$$

for the stochastic clock. This is the key idea for obtaining the representation

$$f_{GY,a}(w) = E^Q\!\left[\frac{e^{2(\nu+1)[w-\tau_{\nu,a}]^+} - 1}{2(\nu+1)} \cdot (R_a^{(\nu)})^2\right],$$

for all $w > 0$. Indeed, fix any positive real number $x$, and consider the process $A^{(\nu)}$ at $x$ on the set of all events where the passage time $\tau_{\nu,a}$ takes values less than or equal to $x$. Break the integral defing $A^{(\nu)}(x)$ at $\tau_{\nu,a}$. The first summand then is $A^{(\nu)}$ at time $\tau_{\nu,a}$ and so is equal to $a$. In the second summand you want to restart the Brownian motion in the exponent of the integrand at $\tau_{\nu,a}$. Thus shift the variable of integration accordingly. The second integral then is the product of $\exp(2 \cdot (B(\tau_{\nu,a}) + \nu\tau_{\nu,a}))$ times $A^{(\nu)}$ at $x - \tau_{\nu,a}$, by abuse of language after having applied Strong Markov. This last process is such that it is independent of the information at time $\tau_{\nu,a}$. Unravelling the definition of $\tau_{\nu,a}$, the first factor is the square of the Bessel process $R^{(\nu)}$ at time $a$. Now taking the expectation conditional on the information at $\tau_{\nu,a}$, we thus get:

$$E^Q\Big[\big(A^{(\nu)}(x)-a\big)^+\,\Big|\,\mathscr{F}_{\tau_{\nu,a}}\Big] = \big(R^{(\nu)}(a)\big)^2 \cdot E^Q\Big[A^{(\nu)}\big(\big[x-\tau_{\nu,a}\big]^+\big)\Big].$$

The $Q$–expectation of $A^{(\nu)}(w)$ is $(\exp(2(\nu+1)w)-1)/(2(\nu+1))$, as shown in §4 or using [**Y**, §4]. On substitution the required expression for $f_{GY,a}$ follows.

At first sight this appears to complicate the problem. However, you will see, it is just what is especially suited to the Laplace transform $F_{GY,a}$ of $f_{GY,a}$:

$$F_{GY,a}(z) = \int_0^\infty e^{-zw} E^Q\left[\frac{e^{2(\nu+1)[w-\tau_{\nu,a}]^+}-1}{2(\nu+1)} \cdot (R_a^{(\nu)})^2\right]dw\,.$$

Still, for computing this integral you want to interchange the Laplace integral with the expectation $E^Q$. If $z$ is real it seems best to follow Yor's proposal for justifying this. Indeed with the integrand of the double integral in question positive and measurable, apply Tonelli's theorem now but justify only in a later step that any of the resulting integrals is finite. The case of a general argument $z$ is reduced to this case considering the absolute value of the integrand, and the result is the identity

$$F_{GY,a}(z) = \frac{1}{z(z-2(\nu+1))} E^Q\Big[e^{-z\tau_{\nu,a}}(R^{(\nu)})^2\Big]$$

of measurable functions for any complex number $z$ with $\mathrm{Re}\,(z) > 2(\nu+1)$. To identify the expectation in the numerator as $D_\nu(a,z)$ write it out as

$$\int_0^\infty x^2 \cdot E^Q\Big[e^{-z\tau_{\nu,a}}\Big|R_a^{(\nu)}=x\Big]\cdot p_{\nu,a}(1,x)\,dx\,,$$

where $p_{\nu,a}$ is the Bessel semigroup of index $\nu$ starting at 1 at time $a$. Under the hypothesis $\nu \geq 0$, Yor has in [**Y80**, Théorème 4.7, p.80] (see also [**GY**, Lemma 2.1 and Proposition 2.6]) explicitly computed the conditional expectation factor of this integrand as the following quotient of modified Bessel functions:

$$E^Q\left[e^{-z\cdot\tau_{\nu,a}}\Big|R^{(\nu)}(a)=w\right] = \frac{I_{\sqrt{2z+\nu^2}}}{I_\nu}\left(\frac{w}{a}\right).$$

Explicit expressions for the Bessel semigroups $p_{\nu,a}$ are known only for $\nu > -1$, see [**Y80**, (4.3), p.78] or [**GY**, Proposition 2.2]. Indeed, in this case the density $p_{\nu,a}(1,w)$ of the Bessel semigroup with index $\nu$ and starting point 1 at time $a$ is

$$p_{\nu,a}(1,w) = \frac{w^{\nu+1}}{a}\cdot e^{-\frac{1+w^2}{2a}}\cdot I_\nu\left(\frac{w}{a}\right).$$

So these results are available both only if $\nu \geq 0$, and this is where the restriction on $\nu$ comes from. On substitution we then get as required

$$D_\nu(a,z) = \frac{e^{-\frac{1}{2a}}}{a}\int_0^\infty e^{-\frac{x^2}{2a}}x^{\nu+3}I_{\sqrt{2z+\nu^2}}\Big(\frac{x}{a}\Big)dx\,.$$

To complete the Tonelli argument proposed to us by Yor and to complete the proof, we have to establish the finiteness of this integral for any fixed complex number $z$ with $\mathrm{Re}\,(z) > 2(\nu+1)$. There is a further technical point to be taken care of. Choose the principal branch of the logarithm to define the square root on the complex plane with the non–positive real line deleted, whence

$$\mu = \sqrt{2z + \nu^2}$$

has a positive real part. Integrability now follows using the asymptotic behaviour of the Bessel function $I_\mu$ near the origin and towards infinity in an essential way. Indeed, from [**L**, §5.7] recall that $I_\mu$ is a continuous function on the positive real line in particular whose asymptotic behaviour near zero is

$$I_\mu(\xi) \approx \frac{\xi^\mu}{2^\mu \Gamma(1+\mu)} \qquad \text{for } \xi \downarrow 0\,.$$

Thus, with $\nu + 3$ and $\mathrm{Re}\,(\mu)$ positive, no integrability problems arise for $x$ near the origin. From [**L**, §5.11] we have for large arguments

$$I_\mu(\xi) \approx \frac{e^\xi}{\sqrt{2\pi\xi}} \qquad \text{for } \xi \to \infty\,,$$

whence the factor $\exp(-x^2/(2a))$ dominates the asymptotic behaviour of the integrand of $D_\nu$ with $x$ to infinity. Thus the integral is finite, completing the computation of the modified Geman–Yor Laplace transform of the Lemma.

# Part III   Complex analytic part of the proof

## 10.  A Lemma about complex square roots

Lifting the non–negativity restriction on $\nu$ of §9 Lemma to the validity of §7 Theorem is by analytic continuation and is achieved in §15. As a prerequisite this section establishes the following result about the complex square root associated to the principal branch of the logarithm

**Lemma:**   *For any complex number $z$ with $\mathrm{Re}\,(z) \geq 2$, we have $\mathrm{Re}\,(\sqrt{2z + \nu^2}) > |\mathrm{Re}\,(\nu)|$ for any complex number $\nu$ with $|\mathrm{Im}\,(\nu)| \leq 1$.*

For the proof first write out the complex numbers involved by explicitly introducing their real and imaginary parts. So let $z = z_0 + iz_1$ for real numbers $z_0 \geq 2$ and $z_1$, and $\nu = \nu_0 + i\nu_1$ with real numbers $|\nu_1| \leq 1$ and $\nu_0$. Then we have

$$2z + \nu^2 = \xi_0 + i2\xi_1 = \xi$$

setting $\xi_0 = 2z_0 + \nu_0^2 - \nu_1^2$ and $\xi_1 = z_1 + \nu_0\nu_1$. Observe that $\xi_0 \geq 3 + \nu_0^2$.

Moreover recall from §16.1 the working of the complex square root. Write $\xi$ as the product of the unique point on the unit circle with angle $|\theta| < \pi$ times the stretching factor $|\xi|$. Then $\xi^{1/2} = |\xi|^{1/2} \exp(i\theta/2)$ whence $\mathrm{Re}\,(\xi^{1/2}) = |\xi|^{1/2} \cos(\theta/2)$. Moreover notice that with $\xi_0$ positive the absolute value of $\theta$ is smaller than 90 degrees.

Distinguish if the absolute value of the argument $\theta$ of $\xi$ is smaller than 60 degress or not, i.e., if $\xi$ is inside the wedge in the right half–plane between the line $y = 2x$ and $y = -2x$ or not.

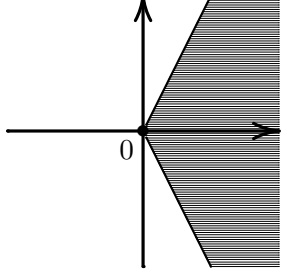

**Figure 2** The 60°–wedge

First let $\xi$ be not inside the shaded 60°–wedge of Figure 2, i.e., let $\xi_0 \leq |\xi_1|$. That the real part of $\xi^{1/2}$ is bigger than the absolute value of real part of $\nu$ then follows by a series of estimates. With the absolute value of $\theta/2$ smaller than 45 degrees, $\cos(\theta/2)$ is bigger than $\cos(\pi/4) = 2^{-1/2}$. For the real part of the square root of $\xi$ we thus get

$$\mathrm{Re}\,(\sqrt{\xi}) = \sqrt{|\xi|} \cdot \cos\left(\frac{\theta}{2}\right) \geq \frac{\sqrt{2}}{2}(\xi_0^2 + 4\xi_1^2)^{1/4}\,.$$

Using $\xi_0^2 \leq \xi_1^2$, we have $\xi_0^2 + 4\xi_1^2 \geq 5\xi_0^2$, whence

$$\mathrm{Re}\,(\sqrt{\xi}) \geq \frac{\sqrt{2}}{2}5^{1/4}\sqrt{\xi_0}\,.$$

Recalling $\xi_0 > 3 + \nu_0^2$, we have $\xi_0^{1/2} > |\nu_0|$. Since moreover $5^{1/4}$ is bigger than $2^{1/2}$,

$$\mathrm{Re}\,(\sqrt{\xi}) > \sqrt{3 + \nu_0^2} > |\nu_0| = |\mathrm{Re}\,(\nu)|,$$

as was to be shown.

Now let $\xi$ be inside the shaded 60°–wedge of Figure 2, i.e., let $|\xi_1| < \xi_0$. Again we have to show that the real part of $\xi^{1/2}$ is bigger than the absolute value of real part of $\nu$. For this develop the square root at $\xi$ into its Taylor series around $\xi_0$:

$$\sqrt{\xi} = \sqrt{\xi_0} + \frac{2i\xi_1}{2\sqrt{\xi_0}} + \frac{1}{2!}\frac{-1}{4 \cdot \eta_\alpha^{3/2}}(2i\xi_1)^2$$

with $\eta_\alpha = \xi_0 + \alpha 2i\xi_1$ for an $\alpha$ in $[0,1]$. Then $\eta_\alpha$ is also inside the 60°–wedge. Thus the absolute value of the angle $\theta_\alpha$ of $\eta_\alpha$ is less than 60 degrees. Hence the absolute value of $3\theta_\alpha/2$ is less than 90 degrees. Thus $\cos(3\theta_\alpha/2)$ is positive, and we have

$$\mathrm{Re}\left(\frac{1}{\eta_\alpha^{3/2}}\right) = \frac{\mathrm{Re}\,(e^{-i3\theta_\alpha/2})}{|\eta_\alpha|^{3/2}} = \frac{\cos(3\theta_\alpha/2)}{|\eta_\alpha|^{3/2}} > 0\,.$$

Using this non–negativity result in the above Taylor expansion,

$$\mathrm{Re}\,(\sqrt{\xi}) = \sqrt{\xi_0} + \frac{\xi_1^2}{2}\mathrm{Re}\left(\frac{1}{\eta_\alpha^{3/2}}\right) > \sqrt{\xi_0}\,.$$

Using that $\xi_0^{1/2}$ is bigger than the square root of $3 + \nu_0^2$ and thus bigger than $|\nu_0|$, this completes the proof of the Lemma.

## 11. First step of the analytic continuation

As a first step in the analytic continuation argument this section studies the analytic properties of the generalized first Weber integral $D_\nu(a, z)$ of §7 Theorem. On choosing the square root associated to the principal branch of the logarithm recall

$$D_\nu(a, z) = \frac{e^{-\frac{1}{2a}}}{a} \int_0^\infty e^{-\frac{x^2}{2a}} x^{\nu+3} I_{\sqrt{2z+\nu^2}}\Big(\frac{x}{a}\Big) dx$$

for any positive real number $a$ and for any complex numbers $z$ and $\nu$ such that the real part of $\nu + 4 + (2z + \nu^2)^{1/2}$ is positive. Taking $z$ with $\mathrm{Re} > 2(\nu+1)$, we know from §9 Lemma that this integral is finite if $\nu$ is non–negative. The key point is that it stays finite if $\nu$ varies in the band of height two which is symmetric with respect to the real axis. Using the confluent hypergeometric function $\Phi$ discussed in [L, §9.9], the precise result is the following

**Lemma:** *For any complex number $z$ with $\mathrm{Re}\,(z) \geq 2$ and any positive real number $a$, we have*

$$D_\nu(a, z) = \Gamma\Big(\frac{\nu+4+\mu}{2}\Big) \cdot \frac{1}{\Gamma(\mu+1)} \Phi\Big(\frac{\nu+4+\mu}{2}, \mu+1; \frac{1}{2a}\Big) \cdot e^{-\frac{1}{2a}} \cdot (2a)^{\frac{\nu+2-\mu}{2}}$$

*for any complex number $\nu$ such that $|\mathrm{Im}\,(\nu)| \leq 1$ and with $\mu = (2z + \nu^2)^{1/2}$.*

**Corollary:** *For any complex number $z$ with $\mathrm{Re}\,(z) \geq 2$ and any positive real number $a$, sending $\nu$ to $D_\nu(a, z)$ defines an analytic map on the set of all complex numbers $\nu$ with $|\mathrm{Im}\,(\nu)| < 1$.*

**Proof of the Lemma:** For proving the Lemma first establish finiteness of $D_\nu(a, z)$. Using the asymptotic behaviour of the modified Bessel function factor of its integrand as in §9, it will be finite if the real part of $\mu(z) = (2z + \nu^2)^{1/2}$ is bigger than the absolute value of the real part of $\nu + 4$. This is implied by §10 Lemma if the imaginary part of $\nu$ is less than or equal to one and the real part of $z$ is greater than or equal to two.

For computing this integral we modify the discussion in [W, §13.3, pp.393f] of Hankel's generalization of Weber's first integral. The idea is to expand the modified Bessel function in the integrand of

$$I = \int_0^\infty e^{-\frac{x^2}{2a}} x^{\nu+3} I_\mu\Big(\frac{x}{a}\Big) dx$$

into its series of [L, §5.7] and integrate term by term. Using [L, §9.9] this is justified by the absolute convergence of the series for the confluent hypergeometric series which is to result, and we get

$$I = \frac{1}{(2a)^\mu} \sum_{n=0}^\infty \frac{1}{\Gamma(\mu+1+n)} \frac{(2a)^{-2n}}{n!} \int_0^\infty e^{-\frac{x^2}{2a}} x^{\nu+3+\mu+2n}\, dx\,.$$

Changing variables $y = (2a)^{-1}x^2$, compute any $n$–th integral as

$$\int_0^\infty e^{-\frac{x^2}{2a}} x^{\nu+3+\mu+2n}\, dx = \frac{1}{2} \cdot \Gamma\Big(\frac{\nu+\mu+4}{2} + n\Big) \cdot (2a)^{\frac{\nu+\mu+4}{2}+n}.$$

Sorting out the series of the pertinent confluent hypergeometric function we thus get

$$I = \frac{1}{2} \cdot \frac{\Gamma\big((\nu+4+\mu)/2\big)}{\Gamma(\mu+1)} \cdot \Phi\Big(\frac{\nu+4+\mu}{2}, \mu+1; \frac{1}{2a}\Big) \cdot (2a)^{\frac{\nu-\mu+4}{2}}.$$

Multiplying this expression with $\exp(-(2a)^{-1})/a$, the Lemma follows.

**Proof of the Corollary:**   For proving the Corollary, using §10 Lemma notice that $\mu(\nu) = (2z+\nu^2)^{1/2}$ is analytic at any $\nu$ with $|\operatorname{Im}(\nu)| < 1$. Thus the confluent hypergeometric function of the Lemma divided by $\Gamma(\mu(\nu)+1)$ is analytic at any such $\nu$ as composition of analytic functions [**L**, §9.9]. Again using §10 Lemma, the real part of $\nu + \mu(\nu)$ is moreover positive for any $\nu$ under consideration. With the gamma function analytic on the right half–plane, thus $\Gamma((\nu + \mu(\nu) + 4)/2)$ is analytic in $\nu$, too, as the composition of analytic functions. The proof of the Corollary is complete.

## 12. Analyticity properties of the moments of $A^{(\nu)}$

This section considers the first two moments $E^Q[(A_x^{(\nu)})^n]$ of Yor's accumulation processes $A^{(\nu)}$ at any positive real number $x$ as functions $e_{n,x}$ in the complex variable $\nu$, and provides the prototype for the discussion to follow. While explicit expressions for these moments are in [**Y**, §4], representations better suited for analytic continuation are obtained by direct re-computation. For any first moment function $e_{1,x}$ we have

$$e_{1,x}(\nu) = \frac{e^{2x(\nu+1)} - 1}{2(\nu+1)} = x \sum_{n=0}^\infty \frac{(2x(\nu+1))^n}{(n+1)!},$$

originally only if $\nu$ is different from $-1$. However, the series is absolutely convergent for all $\nu$ and converges uniformly for $\nu$ in any compact subset of the complex plane. Thus $e_{1,x}$ is analytic in $\nu$ on the whole complex plane with $e_{1,x}(-1) = x$.

Any second moment function $e_{2,x}$ is given by

$$e_{2,x}(\nu) = \frac{2}{2(\nu+1)} e^{2x(\nu+1)} \frac{e^{2x(\nu+3)} - 1}{2(\nu+3)} - \frac{2}{2(\nu+1)} \frac{e^{4x(\nu+2)} - 1}{4(\nu+2)}$$

for any complex number $\nu$ different from $-1$, $-2$, and $-3$. It is an analytic function on the whole complex plane whose extensions over the bad points are given at $\nu = -1$ by

$$e_{2,x}(\nu) = 2\frac{e^{4x(\nu+2)} - 1}{4(\nu+2)},$$

at $\nu = -2$ by

$$e_{2,x}(\nu) = \frac{2}{2(\nu+1)} e^{2x(\nu+1)} \frac{e^{2x(\nu+3)} - 1}{2(\nu+3)} - \frac{2x}{2(\nu+1)},$$

and finally at $\nu = -3$ by

$$e_{2,x}(\nu) = \frac{2x}{2(\nu+1)} e^{2x(\nu+1)} - \frac{2}{2(\nu+1)} \frac{e^{4x(\nu+2)} - 1}{4(\nu+2)}.$$

## 13. Existence and analyticity properties of the functions $f_{GY,a}$

In this section any of the auxiliary functions $f_{GY,a}$ of §6 is considered as function in the variable $\nu$

$$L(\nu) = E^Q\big[\big(A_x^{(\nu)} - a\big)^+\big]$$

for any fixed positive real numbers $a$ and $x$. Using §9 Lemma, we know it is defined for non–negative real numbers $\nu$. However, this has been achieved in a very indirect way only: for these values of $\nu$ the Laplace transforms of the corresponding functions $f_{GY,a}$ have been shown to be finite. Now more is true indeed

**Lemma:** *The function $L$ extends to a function on the complex plane which is analytic at each point, and for which we have the majorization*

$$|L(\nu)| \le e^{\frac{x}{2}\operatorname{Im}^2(\nu)} \cdot E^Q[A_x^{(\operatorname{Re}(\nu))}].$$

For the proof of the Lemma now set $f(w) = (w-a)^+$. Applying Girsanov such that $W_x = \nu x + B_x$ becomes a standard Brownian motion, and, dropping reference to this new measure, we get

$$L(\nu) = E\Big[f\big(A_x^{(0)}\big)e^{\nu W_x}\Big] \cdot e^{-\frac{x}{2}\nu^2}.$$

For establishing the analyticity statement of the Lemma it is thus sufficient to show that the expectation factor is analytic in any complex number $\nu$. This is true by definition if we have the convergent series

$$E\Big[f\big(A_x^{(0)}\big)e^{\nu W_x}\Big] = \sum_{m=0}^{\infty} \frac{\nu^m}{m!} E\Big[f\big(A_x^{(0)}\big)W_x^m\Big]$$

for all $\nu$. For this it is sufficient to show that the series is absolutely convergent for all $\nu$. Using the Cauchy–Schwarz inequality this is implied by the convergence of

$$\sum_{m=0}^{\infty} \frac{|\nu|^m}{m!} \sqrt{E\Big[f^2\big(A_x^{(0)}\big)\Big]} \sqrt{E\big[W_x^{2m}\big]}$$

for all $\nu$. Herein $E[f^2(A_x^{(0)})]$ is majorized by the second moment of Yor's zero drift accumulation process at $x$, and so is finite. Since $E[W_x^{2m}]$ is majorized by $\pi^{-1/2} \cdot (2x)^m \cdot m!$ for all $m \ge 0$, convergence follows using the ratio test. Actually, we have so established yet another upper bound to the price of the Asian option.

To establish the majorization of the Lemma, taking absolute values inside the expectation in the above Girsanov representation of $L$ gives:

$$|L(\nu)| \le E\Big[f\big(A_x^{(0)}\big) \cdot |e^{\nu W_x}|\Big] \cdot |e^{-\frac{x}{2}\nu^2}|.$$

The absolute value of the exponential factors are the exponentials of the real parts of the respective arguments. Majorizing the function of Yor's zero drift accumulation process by this process itself, we get

$$|L(\nu)| \leq e^{\frac{x}{2}\operatorname{Im}^2(\nu)} \cdot E\left[A_x^{(0)} e^{\operatorname{Re}(\nu)W_x}\right] e^{-\frac{x}{2}\operatorname{Re}^2(\nu)}.$$

Reversing the Girsanov then completes the proof of the Lemma.

## 14. Existence and analyticity properties of the functions $F_{GY,a}$

While the preceding section has studied the expectations

$$L(x, \nu) = E^Q\left[\left(A_x^{(\nu)} - a\right)^+\right]$$

for any fixed positive real numbers $a$ and $x$ as function in the complex variable $\nu$ only, this section moreover treats $x$ as a variable. If $\nu$ is any non–negative real number, we have from §9 Lemma that the integrals of the Laplace transforms

$$F(\nu)(z) = \int_0^\infty e^{-zx} L(x, \nu)\,dx$$

are finite if $\operatorname{Re}(z) > 2(\nu+1)$. In this section we give independent proofs of the following two more general statements

**Lemma:**  *For any complex number $\nu$, the Laplace transform $F(\nu)(z)$ is finite for any complex number $z$ with $\operatorname{Re}(z) > \max\{0, \operatorname{Im}^2(\nu)/2 + 2(\operatorname{Re}(\nu)+1)\}$.*

**Corollary:**  *For any complex number $z$ with a positive real part, the map sending $\nu$ to $F(\nu)(z)$ is analytic in all complex numbers $\nu$ with $\operatorname{Re}(z) > \operatorname{Im}^2(\nu)/2 + 2(\operatorname{Re}(\nu)+1)$.*

**Proof of the Lemma:**  As first step in proving the Lemma, we establish for any complex number $\nu$ the majorization

$$|F(\nu)(z)| \leq \int_0^\infty e^{-\left(\operatorname{Re}(z) - \frac{1}{2}\operatorname{Im}^2(\nu)\right)x} E^Q[A_x^{(\operatorname{Re}(\nu))}]\,dx$$

for any complex number $z$ in the sense of measurable functions. Indeed, majorize the absolute value of $F(\nu)(z)$ by taking the absolute value inside the defining integral. The absolute value of the exponential factor then is equal to $\exp(-\operatorname{Re}(z)x)$. Majorizing the absolute value of $L(x, \nu)$ using §13 Lemma, the estimate follows.

Setting $\nu_0 = \operatorname{Re}(\nu)$, now let $\operatorname{Re}(z)$ be positive and such that $\xi_0 = \operatorname{Re}(z) - \operatorname{Im}(\nu)^2/2$ is bigger than $2(\nu_0+1)$. We compute the Laplace transform of the right hand side of the above inequality using §12. If $\nu_0$ is different from minus one,

$$\int_0^\infty e^{-\xi_0 x} E^Q[A_x^{(\nu_0)}]\,dx = \frac{1}{\xi_0(\xi_0 - 2(\nu_0+1))}$$

using that $\xi_0$ is bigger than $2(\nu_0+1)$ to compute the improper integrals. It converges to $\xi_0^{-2}$ with $\nu_0$ going to minus one. Thus it is seen to coincide with the Laplace transform for the case $\nu_0 = -1$. The Laplace transforms are finite if $\xi_0$ is positive and bigger than $2(\nu_0+1)$, and then $F(\nu)(z)$ is finite *a forteriori*. The proof of the Lemma is complete.

**Proof of the Corollary:**   There are a number of reductions for the proof of the Corollary. First, it is sufficient to show that $F(\nu)(z)$ is one–fold partially complex differentiable with respect to $\nu$. For this we can assume that $\nu$ belongs to any sufficiently small bounded open subset $U$ of the complex plane. We claim that computing this partial derivative then is by differentiation under the integral defining $F(\nu)(z)$. This will complete the proof since the integrand of $F(\nu)(z)$ is analytic in $\nu$ using §13 Lemma. Thus we are finally reduced to show that the restriction of the partial derivative of the integrand of $F(\nu)(z)$ to $U$ times the positive real line can be majorized in absolute value by an integrable function on the positive real line.

As a first step, we establish for the expectation $L(x,\nu)$ the majorization

$$e^{-\frac{x}{2}\operatorname{Im}^2(\nu)} \cdot \left| \frac{\partial L}{\partial \nu}(x,\nu) \right| \leq |\nu| x \cdot E^Q\big[ A_x^{(\operatorname{Re}(\nu))} \big] + \sqrt{E^Q\big[ B_x^2 \big]} \cdot \sqrt{E^Q\big[ \big( A_x^{(\operatorname{Re}(\nu))} \big)^2 \big]} \, .$$

for any positive real $x$ and any $\nu$ in $U$.

Indeed, apply Girsanov such that $W_x = \nu x + B_x$ becomes a standard Brownian motion. Dropping reference to this new measure,

$$L(x,\nu) = E\Big[ f\big( A_x^{(0)} \big) e^{\nu W_x} \Big] \cdot e^{-\frac{x}{2}\nu^2}$$

with $f(w) = (w - a)^+$. The power series in $\nu$ of the expectation factor derived in §13 is uniformly convergent on compact subsets of the complex plane. On $U$ differentiation of this series with respect to $\nu$ so is term by term yielding

$$\frac{\partial}{\partial \nu} E\Big[ f\big( A_x^{(0)} \big) e^{\nu W_x} \Big] = E\Big[ f\big( A_x^{(0)} \big) W_x e^{\nu W_x} \Big] \, .$$

The partial derivative of $L$ with respect to $\nu$ thus equals

$$\frac{\partial L}{\partial \nu}(x,\nu) = E\Big[ f\big( A_x^{(0)} \big) W_x e^{\nu W_x} \Big] \cdot e^{-\frac{x}{2}\nu^2} - x\nu \cdot L(x,\nu)$$

on reversing the Girsanov in the second summand. To estimate the absolute value of the second summand use the majorization of §13 Lemma. To estimate the absolute value of the expectation in the first summand, apply Cauchy–Schwarz to get

$$\left| E\Big[ f\big( A_x^{(0)} \big) W_x e^{\nu W_x - \frac{x}{2}\nu^2} \Big] \right|^2 \leq E\big[ W_x^2 \big] \cdot E\Big[ \Big| f\big( A_x^{(0)} \big) e^{\nu W_x - \frac{x}{2}\nu^2} \Big|^2 \Big] \, .$$

Majorize the second expectation on the right hand side arguing as for §13 Lemma, and the desired inequality follows.

As a next step, we show that after shrinking $U$ if necessary the absolute values of the functions $x \mapsto \exp(-zx)(\partial L / \partial u)(x,u)$ with $u$ in $U$ can be majorized in dependency on $U$ by a function which is integrable on the positive real line.

For this apply the inequality of the first step. There we have $E^Q\big[B_x^2\big] = x$, for any positive real number $x$. Using §12 for the other expectations, for any $u$ in $U$, the absolute value of the partial derivatives $x \mapsto (\partial L/\partial u)(x, u)$ can thus be majorized by a function $f_u$ on the positive real line with the following properties. It is continuous in both the parameter $\mathrm{Re}\,(u)$ and its argument $x$, it behaves like a polynomial for small values of $x$, and for large values of $x$ behaves like $\exp(m_u x)$ on abbreviating $m_u = \max\{0, \mathrm{Im}^2(u)/2 + 2(\mathrm{Re}\,(u)+2)\}$. Thus $x \mapsto \exp(-zx) \cdot (\partial L/\partial u)(x, u)$ is integrable on the positive real line if $\mathrm{Re}\,(z)$ is positive and bigger than $\mathrm{Im}^2(u)/2 + 2(\mathrm{Re}\,(u)+2)$.

Now recall we had fixed any $\nu$ such that $\mathrm{Re}\,(z) > \mathrm{Im}^2(\nu)/2 + 2(\mathrm{Re}\,(\nu)+2)$ and then chosen a bounded open neighborhood $U$ of $\nu$. Shrink $U$ if necessary so that $\mathrm{Re}\,(z) > \mathrm{Im}^2(u)/2 + 2(\mathrm{Re}\,(u)+2)$ for any $u$ in its closure. Choose $u^*$ in this closure with $m_{u^*}$ maximal. This takes care of the growth behaviour of the exponential factors in the $f_u$ in such a way that the ones in $f_{u^*}$ dominate those of all other functions $f_u$ for larger values of $x$. Now take stock of the functions that occur in $f_{u^*}$ as coefficients of the exponentials. They are continuous in $\mathrm{Re}\,(u)$ in particular. Letting $u$ vary over the compact closure of $U$ majorize them by constants also depending on $U$ so that the resulting function $f_U$ majorizes all other functions $f_u$ with $u$ in $U$. Since $f_U$ behaves for large values of the argument $x$ like $\exp(m_{u^*} x)$ too, $x \mapsto \exp(-x\mathrm{Re}\,(z)) \cdot f_U(x)$ is integrable on the postive real line by the choice of $u^*$. So it qualifies as a majorizing function we have been looking for, and this completes the proof of the Corollary.

## 15. A summing up

It remains to pull things together and establish the two results of §7.

First notice that §7 Corollary is implied by §7 Theorem using the computation of $D_\nu(a, z)$ in §11 Lemma. Thus we are reduced to prove §7 Theorem.

The proof of §7 Theorem is by analytic continuation using the identity theorem recalled in §16.2. As a consequence of §9 Lemma it remains to establish this result for negative indices $\nu$ only. The existence of the Laplace transform on all complex numbers $z$ with positive real part required in the Theorem then has been proved in §14 Lemma.

For establishing the crucial Laplace transform identity of §7 Theorem, thus let $z$ be any complex number with real part $\mathrm{Re}\,(z) > 4$ and choose $0 < \varepsilon < 1$. Using §11 Corollary, the generalized first Weber integral

$$D_\nu(a, z) = \frac{e^{-\frac{1}{2a}}}{a} \int_0^\infty e^{-\frac{x^2}{2a}} x^{\nu+3} I_{\sqrt{2z+\nu^2}}\Big(\frac{x}{a}\Big) dx$$

is analytic in $\nu$ on the $\varepsilon$–thickened real line $A_\varepsilon$ which consists of all complex numbers $\nu$ with $|\mathrm{Im}\,(\nu)| < \varepsilon$. Picture $A_\varepsilon$ as the band of height $2\varepsilon$ symmetric with respect to the real axis. If we choose $\varepsilon$ so small that $\mathrm{Re}\,(z) > 2(2+2\varepsilon)$, we claim to have analyticity of the Laplace transform

$$F(\nu)(z) = \int_0^\infty e^{-zx} L(x, \nu) \, dx$$

as a function in $\nu$ on the subset $B_\varepsilon$ of $A_\varepsilon$ which consists of all complex numbers $\nu$ with $|\text{Im}(\nu)| < \varepsilon$ and $\text{Re}(\nu) < \varepsilon$. Indeed, if $\text{Re}(\nu) < \varepsilon$, we have $2(2+2\varepsilon) > 2\varepsilon+2(\text{Re}(\nu)+2)$. If $|\text{Im}(\nu)| < \varepsilon$, we have $2\varepsilon > \text{Im}^2(\nu)/2$ since $\varepsilon < 1$. Any $\nu$ in $B_\varepsilon$ so satisfies the inequality $\text{Re}(z) > \text{Im}^2(\nu)/2 + 2(\text{Re}(\nu)+2)$ of §14 Corollary, and the claim follows.

So the functions $F(\nu)(z)$ and

$$D_\nu^*(a,z) = \frac{D_\nu(a,z)}{z(z-2(\nu+1))}$$

are analytic as functions in $\nu$ on the $\varepsilon$–thickened half–line $B_\varepsilon$. It is now a consequence of §10 Lemma that we have $F(\nu)(z) = D_\nu^*(a,z)$ for all $\nu \geq 0$ such that $2(\nu+1) < \text{Re}(z)$. With $\text{Re}(z) > 4$ this so holds *a fortiori* for all $\nu \geq 0$ such that $2(\nu+1) < 4$, i.e., for all non–negative real numbers $\nu$ smaller than 1. With $\varepsilon$ smaller than 1 we so have $F(\nu)(z) = D_\nu^*(a,z)$ for any $\nu$ in the subinterval $(0,\varepsilon)$ of $B_\varepsilon$. With $B_\varepsilon$ open and connected the identity theorem literally applies to give that $F(\nu)(z) = D_\nu^*(a,z)$ for all $\nu$ in $B_\varepsilon$. This identity then holds *a fortiori* for all real numbers $\nu$ in $B_\varepsilon$, i.e., for all $\nu < \varepsilon$. The proof of §7 Theorem is complete.

# 16. Appendix: notions from complex analysis

This section reviews pertinent results about complex analysis. It is not meant as a substitute for a closer study of basic texts like [**FB**], [**C**], or [**R**].

On choosing any root $i$ of $x^2+1$, any complex number $z$ has the form $z = x+iy$ where the real number $x = \text{Re}(z)$ is its real and the real number $y = \text{Im}(z)$ its imaginary part. The complex conjugate of $z$ has the same real part as $z$, its imaginary part is the negative of that of $z$. The product of $z$ and its complex conjugate is a non–negative real number whose square root is the absolute value $|z|$ of $z$. Thinking of $z$ as the pair with entries the real and imaginary part, it equals the Euclidian distance of $z$ to the origin in the plane.

## 16.1. Exponential function, logarithms, and square roots

In analogy to analysis over the reals we have the notion of a complex–valued function $f$ being differentiable in an element of some open subset of the complex plane $\mathbf{C}$ where it is defined. In contrast to real analysis, however, complex differentiability is equivalent to *analyticity* in the sense that if $f$ is differentiable in a point of its domain of definition, it is already differentiable to any order there and the resulting Taylor series converges in a neighborhood of that point. A basic example of an analytic function is the *exponential function* exp given by the series $\exp(z) = \sum_{n=0}^{\infty} z^n/n!$ for any complex number $z$. It extends the exponential function on the real line. To get accustomed to its properties, the absolute value of the exponential function at any complex number $z$ equals $\exp(\text{Re}(z))$, the value of the exponential function at the real part of $z$:

$$|e^z| = e^{\text{Re}(z)}.$$

Indeed, the square of the absolute value of $\exp(z)$ is the product of $\exp(z)$ with its complex conjugate. Using the series, the complex conjugate of $\exp(z)$ is the exponential function at the complex conjugate of $z$, and the above relation follows. The fundamental problem now is that this function is periodic and so has no inverse. Indeed, performing the series calculations,

$$\exp(i\theta) = \cos(\theta) + i\sin(\theta),$$

for any angle $\theta$. There is a way out of this. Think of $\exp(i\theta)$ as the element on the unit circle corresponding to the angle $\theta$. Any complex number $z$ different from zero so is a product of its distance $|z|$ from the origin times an element $\exp(i\theta)$, just draw a ray through it and the origin to see this. So $z = \exp(i\theta + \log|z|)$ which forces

$$\log(z) = i\theta + \log|z|.$$

The problem is that $\theta$ is determined only up to integer multiples of $2\pi$. To get rid of this indeterminancy put restrictions on the range of values of the angles. In the sequel the convention $|\theta| < \pi$ is chosen. The upshot is that so you get a continuous logarithm not on the whole complex plane but only on $\mathbf{C}_- = \mathbf{C} \setminus \mathbf{R}_{\leq 0}$, the complex plane cut along non–positive real line. This is the *principal branch* of the logarithm. For any element $z$ there is now a unique angle $\theta = \arg(z)$ of absolute value less than $\pi$, the *argument* of $z$, such that $z = |z|\exp(i\theta)$. This yields the square root:

$$\sqrt{z} = \sqrt{|z|} \cdot e^{i\frac{\arg(z)}{2}}.$$

The effect of this square root associated to the principal branch of the logarithm can be pictured as follows. Any element in $\mathbf{C}_-$ is on a unique ray through the origin with angle $\theta$ between minus $\pi$ and $\pi$ in a distance $R$ from the origin. The effect of taking its square root is then halving the angle of this ray, so putting things into the right half plane, and changing on this new ray the distance of the image under the square root to the origin to the square root of the old distance $R$ [**FB**, I§2].

In particular notice what happens at any point of the cut. The limit from the upper edge of the cut to this point of the square root differs from that taken from the lower edge of the cut by a factor minus one. Indeed, the angles corresponding to the upper and the lower edge of the cut differ by $2\pi$ which upon taking square roots produces the factor $\exp(i\pi)$ and thus gives the minus one by which the limits differ.

### 16.2. The identity theorem

The *identity theorem* asserts that *any two functions analytic on an open and connected subset $G$ of $\mathbf{C}$ are identical if they agree on a sequence in $G$ which converges in $G$.* Indeed, it is sufficient to prove that any function $f$ as in the theorem is identical to zero if it vanishes on any convergent sequence $d_n$ as in the theorem. Suppressing some purely topological nicities, the connectedness of $G$ implies that it is sufficient to prove $f$ identical zero in a on a neighborhood of the limit $d$ of this sequence. So develop $f$ in its Taylor series around $d$. This Taylor series converges on an open neighborhood $B$ of $d$ in $G$. Since $f$ is continuous and vanishes on the points $d_n$, it is zero at $d$ whence the zeroth above Taylor coefficient is zero. Dividing $f$ by $z - d$ so gives a function on $B$ which is analytic there. It equals zero on all $d_n$ different from $d$, and so is zero on $d$, too. Consequently the first Taylor coefficient of $f$ is zero. Proceeding inductively, all coefficients in the above Taylor expansion of $f$ are zero. Thus $f$ is zero on $B$.

## 17. Appendix: Laplace transform

This section collects a number of relevant facts about the Laplace transform from [**B**], [**D**]. The class of functions considered in the sequel is that of functions of *exponential type*, i.e., of continuous, real–valued functions $f$ on the non–negative real line such that there is a real number $a$ for which $\exp(at)f(t)$ is bounded for any $t > 0$. The *Laplace transform* is the linear operator $\mathscr{L}$ that associates to any function $f$ of exponential type the complex–valued function $\mathscr{L}(f)$, analytic on a suitable complex half–plane, which for any complex number $z$ with sufficiently big real part is explicitly given by:

$$\mathscr{L}(f)(z) = \int_0^\infty e^{-zt} f(t)\, dt.$$

To fix ideas, consider the function given by $f(t) = \exp(at)$ with $a$ any positive real number so that formally $\mathscr{L}(f)(z) = \int_0^\infty \exp(-(z-a)t)\, dt$. Recall from §16.1 that the absolute value of $\exp(w)$ is $\exp(\operatorname{Re}(w))$, for any complex number $w$. This improper integral thus exists iff the real part of $z - a$ is positive, i.e., iff the real part of $z$ is bigger than $a$. In this case, $\mathscr{L}(f)(z)$ equals $(z-a)^{-1}$ and is analytic on the half plane $\{z | \operatorname{Re}(z) > a\}$. Now mind the following fallacy! The function $z \mapsto (z-a)^{-1}$ is defined for any complex number $z$ different from $a$. However, this function has no relations whatsoever with the Laplace transform on this half–plane. In fact, with the defining integral of the Laplace transform being infinity on $\{z | \operatorname{Re}(z) < a\}$, the latter does not make sense there at all.

That $\mathscr{L}$ is an injection is a consequence of ultimately the Weierstraß approximation theorem. Its inverse, the inverse Laplace transform $\mathscr{L}^{-1}$, is expressed as a contour integral by the *complex inversion formula* of Riemann. It applies to any function $h$ analytic on the half–plane $\{\operatorname{Re}(z) \geq z_0\}$ with $z_0$ any sufficiently big positive real number. Then it asserts that for any positive real number $t$ we have:

$$\mathscr{L}^{-1}(h)(t) = \frac{1}{2\pi i} \int_{z_0 - i\infty}^{z_0 + i\infty} e^{zt} h(z)\, dz\,,$$

if $h$ moreover satisfies a growth condition at infinity such that the above integral exists. For any function $f$ of exponential type we so have in particular $(\mathscr{L}^{-1} \circ \mathscr{L})(f) = f$.

Peter Carr
Banc of America Securities, LLC
9 West 57th Street, 40th floor
New York, NY 10019
USA

*E-mail address*: `pcarr@bofasecurities.com`

Michael Schröder
Lehrstuhl Mathematik III
Universität Mannheim
D–68131 Mannheim
Germany

*E-mail address*: `schroeder@math.uni-mannheim.de`